\documentclass{amsart}
\usepackage{array}
\usepackage{ifthen}
\usepackage{url}
\usepackage[margin=3cm]{geometry}  
\usepackage{amssymb,amsmath, amsfonts}
\newif\iffull

\newtheorem{theorem}{Theorem}[section]
\newtheorem{lemma}[theorem]{Lemma}

\theoremstyle{definition}

\theoremstyle{remark}

\numberwithin{equation}{section}

\def\QQ{\mathbb Q}
\def\ZZ{\mathbb Z}

\newcommand{\klammern}[4][]%
{\ifthenelse{\equal{#1}{}}{\left#2}{\csname#1\endcsname#2}%
#4\ifthenelse{\equal{#1}{}}{\right#3}{\csname#1\endcsname#3}}
\newcommand{\betrag}[2][]{\klammern[#1]{\lvert}{\rvert}{#2}}

\newcommand{\conj}[1]{^{(#1)}}

\newcommand{\abs}[1]{\left\vert#1\right\vert}

\renewcommand{\epsilon}{\DONOTUSE}  %% use \varepsilon instead

\renewcommand{\vartheta}{\DONOTUSE} %% use \theta instead.
\date{\today}

\begin{document}

\title[Difference between powers of two consecutive Balancing numbers]{\small An exponential Diophantine equation related to the difference between powers of two consecutive Balancing numbers}

\author{Salah Eddine Rihane}
\address{Universit\'e des Sciences et de la Technologie Houari-Boumedi\`{e}ne (USTHB), Facult\'e de Math\'ematiques, Laboratoire d'Alg\`{e}bre et Th\'eorie des Nombres, BP 32, 16111 Bab-Ezzouar Alger, Alg\'erie}
\email{salahrihane@hotmail.fr}

\author{Bernadette Faye}
\address{Department of Mathematics, University Gaston Berger of Saint-Louis (UBG), Saint-Louis, Senegal}
\email{bernadette@aims-senegal.org}

\author{Florian Luca}
\address{School of Mathematics, University of the Witwatersrand, Private Bag 3, Wits 2050,
Johannesburg, South Africa; Max Planck Institute for Mathematics,  Vivatsgasse 7, 53111 Bonn, Germany;  
Department of Mathematics,   Faculty of Sciences, University of Ostrava, 30 Dubna 22, 701 03
Ostrava 1, Czech Republic}
\email{Florian.Luca@wits.ac.za}

\author{Alain Togb\'e}
\address{Department of Mathematics, Statistics, and Computer Science, Purdue University Northwest, 1401 S, U.S. 421, Westville IN 46391, USA}
\email{ atogbe@pnw.edu}

%\author{Alain Togb\'e}
%\address{Department of Mathematics, Statistics, and Computer Science Purdue University Northwest, 1401 S, U.S. 421, Westville IN 46391 USA}
%\email{atogbe@pnw.edu}

\thanks{B. Faye was partially supported by a grant from the Simons Foundation}
\thanks{A. Togb\'e was supported in part by Purdue University Northwest.}

%    General info
\subjclass[2010]{11B39, 11J86.}
%\date{Received by the editor March 3, 1998 and, in revised form, April 28, 1998.}
%\dedicatory{This paper is dedicated to our authors.}

\keywords{Pell numbers, Linear form in logarithms, reduction method.}

%\typeset{\LATEX\}

\begin{abstract}
In this paper, we find all solutions of the exponential Diophantine equation $B_{n+1}^x-B_n^x=B_m$ in positive integer variables $(m, n, x)$, where $B_k$ is the $k$-th term of the Balancing sequence.
\end{abstract}

\maketitle

%&&&&&&&&&&&&&&&&&&&&&&&&&%
\section{Introduction}\label{sec1}
%&&&&&&&&&&&&&&&&&&&&&&&&&%
The first definition of balancing numbers is essentially due to Finkelstein \cite{Finkelstein:1965}, although he called them numerical centers. A positive integer $n$ is called balancing number if
$$
1+2+\ldots+(n-1)=(n+1)+(n+2)+\ldots+(n+r)
$$
holds for some positive integer $r$. Then $r$ is called \textit{balancer} corresponding to the balancing number $n$.  For example $6$ and $35$ are balancing numbers with balancers $2$ and $14$, respectively. The $n$-th term of the sequence of balancing numbers is denoted by $B_n$. The balancing numbers satisfy the recurrence relation
$$
B_{n}=6B_{n-1}-B_{n-2},
$$
where the initial conditions are $B_0 = 0$ and $B_1 = 1$. Its first terms are 
$$0, 1, 6, 35, 204, 1189, 6930, 40391, 235416, 1372105, \ldots$$
\medskip
It is well-known that
$$
B_{n+1}^2-B_n^2 =B_{2n+2}, \text{ for any } n\geq 0.
$$
In particular, this identity tells us that the  difference between the square of two consecutive Balancing numbers is still a Balancing number. So, one can ask if this identity can be generalized? 
\medskip

Diophantine equations involving  sum or difference of powers of two consecutive members of a given linear recurrent sequence $\{U_n\}_{n\geq 1}$ were also considered in several papers. For example, in \cite{Marques-Togbe:2010}, Marques and Togb\'e proved that if $s\geq 1$ an integer such that $F_m^s+ F_{m+1}^s$ is a Fibonacci number for all sufficiently large $m$, then $s\in\{1,2\}$. In \cite{Luca-Oyono:2011}, Luca and Oyono proved that there is no integer $s\geq 3$ such that the sum of $s$th powers of two consecutive Fibonacci numbers is a Fibonacci number.  Later, their result have been extended in \cite{Ruiz-Luca:2000} to the generalized Fibonacci numbers and recently in \cite{RFLT:2018} to the Pell sequence.
\medskip

Here, we apply the same argument as in \cite{Luca-Oyono:2011} to the Balancing sequence and prove the following:

\begin{theorem}\label{thm:principal}
The only nonnegative integer solutions  $(m,n,x)$ of the Diophantine equation 
\begin{equation}\label{eq:principal}
B_{n+1}^x-B_n^x=B_m
\end{equation} 
are $(m,n,x)=(2n+2,n,2), (1,0,x), (0,n,0)$.
\end{theorem}

Our proof of Theorem 1.1 is mainly based on linear forms in logarithms of algebraic numbers and a reduction algorithm originally introduced by Baker and Davenport in \cite{Baker-Davenport:1969}. Here, we will use a version due to Dujella and Peth\H{o} in \cite[Lemma 5(a)]{Dujella-Petho:1998}.

%For the proof of Theorem 1.2, we refer the reader to the proof of [\cite{RFLT:2018}, Theorem 1.2], since the argument of the proofs are similar.
%&&&&&&&&&&&&&&&&&&&&&&&&&%
\section{Preliminary results}\label{sec2}
%&&&&&&&&&&&&&&&&&&&&&&&&&%

%%%%%%%%%%%%%%%%%%
\subsection{Balancing sequences} 
%%%%%%%%%%%%%%%%%%
Let $(\alpha,\beta)=(3+2\sqrt{2},3-2\sqrt{2})$ be the roots of the characteristic equation $x^2-6x+1 = 0$ of the Balancing sequence $(B_n)_{n\geq 0}$. The Binet formula for $B_n$ is
\begin{equation}\label{eq:Binet-Bal}
B_n=\dfrac{\alpha^n-\beta^n}{4\sqrt{2}} \quad \text{for all } n\geq 0.
\end{equation}
This implies that the inequality 
\begin{equation}\label{eq:enca-Bal}
\alpha^{n-2} \leq B_n \leq \alpha^{n-1}
\end{equation}
holds for all positive integers $n.$  It is easy to prove that
\begin{equation}\label{eq:upper-bound_for-Bm/Bm+1}
\dfrac{B_n}{B_{n+1}} \leq \dfrac{5}{29}
\end{equation}
holds for any $n\geq 2$.

%%%%%%%%%%%%%%%%%%%%
\subsection{Linear forms in logarithms}
%%%%%%%%%%%%%%%%%%%%

For any non-zero algebraic number $\gamma$ of degree $d$ over $\QQ$, whose minimal polynomial
over $\ZZ$ is $a\prod_{i=1}^d \left(X-\gamma\conj i \right)$, we denote by
\[
h(\gamma) = \frac{1}{d} \left( \log|a| + \sum_{i=1}^d \log\max\left(1,
\betrag{\gamma\conj i}\right)\right)
\]
the usual absolute logarithmic height of $\gamma$.

With this notation, Matveev proved the following theorem (see \cite{Matveev:2000}).

\begin{theorem}\label{lem:Matveev}
Let $\gamma_1,\ldots ,\gamma_s$ be a real algebraic numbers and let  $b_1,\ldots,b_s$ be nonzero rational integer numbers. Let $D$ be the degree of the number field $\QQ(\gamma_1,\ldots ,\gamma_s)$ over $\QQ$ and let $A_j$ be a positive real number satisfying  
$$
A_j=\max\{Dh(\gamma_j),|\log\gamma_j|,0.16\}\quad\text{for}\quad j=1,\ldots ,s. 
$$
Assume that
$$
B\geq\max\{|b_1|,\ldots,|b_s|\}.
$$
If $\gamma_1^{b_1}\cdots\gamma_s^{b_s}-1\neq0$, then
$$
|\gamma_1^{b_1}\cdots\gamma_s^{b_s}-1|\geq\exp(-1.4\cdot 30^{s+3}\cdot s^{4.5}\cdot D^2(1+\log D)(1+\log B) A_1\cdots A_s).
$$
\end{theorem} 
\medskip

%%%%%%%%%%%%%%%%%%%
\subsection{Reduction algorithm}
%%%%%%%%%%%%%%%%%%%

\begin{lemma}\label{le:Baker-Davenport}
Let $M$  be a positive integer, let $p/q$ be a convergent of the continued fraction of the irrational $\gamma$ such that $q > 6M$, and let $A,B,\mu$ be some real numbers with $A > 0$ and $B > 1$. Let
$$
\varepsilon=||\mu q||-M\cdot||\gamma q||,
$$
where $||\cdot||$ denotes the distance from the nearest integer. If    $\varepsilon>0$, then there is no solution of the inequality
$$
0 < m \gamma - n + \mu <A B^{-k}
$$
in positive integers $m,n$ and $k$ with 
$$
m\leq M \quad \text{and} \quad k\geq\dfrac{\log(Aq/\varepsilon)}{\log B}.
$$ 
\end{lemma}

%%%%%%%%%%%%%%%%%%%%%%%%%%%%%%%%%%%%%
\section{The proof of Theorem \ref{thm:principal} }\label{sec3}
%%%%%%%%%%%%%%%%%%%%%%%%%%%%%%%%%%%%%

%%%%%%%%%%%%%%%%%%%%%%%%%%%%%%%%%%%%%
\subsection{An inequality for $x$ versus $m$ and $n$}\label{subsec1}
%%%%%%%%%%%%%%%%%%%%%%%%%%%%%%%%%%%%%
The case $nx=0$ is trivial so we assume that $n\ge 1$ and that $x\ge 1$. Observe that since $B_n< B_{n+1}-B_n < B_{n+1}$, the Diophantine equation \eqref{eq:principal} has no solution when $x=1$. 

When $n=1$, we get $B_m=6^x-1.$ In this case, we have that $m$ is odd. Thus, using the Binet formula \eqref{eq:Binet-Bal}, we obtained the following factorization 
\begin{equation}
\label{eq:except1}
6^x=B_m+1=B_m+B_1=B_{(m+1)/2}C_{(m-1)/2},
\end{equation}
where $\{C_m\}_{m\geq 1}$ is the Lucas Balancing sequence given by the recurrence  $C_{m}=6C_{m-1}-C_{m-2}$ with initial conditions $C_0=2$, $C_1=6$. The Binet formula of the Lucas Balancing sequence is given by $C_n=\alpha^n+\beta^n$. This shows that the largest prime factor of $B_{(m+1)/2}$ is $3$ and by Carnichael Primitive Divisor Theorem we conclude that $(m+1)/2\le 12$, so $m\le 23$. One now checks 
all such $m$ and gets no additional solution with $n=1$. 

\medskip

So, we can assume that $n\geq 2$ and $x\geq 3$. Therefore, we have
$$
B_m=B_{n+1}^x-B_n^x\ge B_3^3-B_1^3= 215,
$$
which implies that $m> 4$. Here, we use the same argument from \cite{Luca-Oyono:2011} to bound $x$ in terms of $m$ and $n$. Since most of the details are similar, we only sketch the argument.
\medskip

Using inequality \eqref{eq:enca-Bal}, we get
$$
\alpha^{m-1}>B_m=B_{n+1}^x-B_n^x\geq B_n^x>\alpha^{(n-2)x},
$$
and
$$
\alpha^{m-2}<B_m= B_{n+1}^x-B_n^x<B_{n+1}^x<\alpha^{nx}.
$$
Thus,
\begin{equation}\label{eq:encad-m-terms-n-x}
(n-2)x+1<m<nx+2.
\end{equation}
Estimate \eqref{eq:encad-m-terms-n-x} is essential for our purpose.

Now, we rewrite the equation \eqref{eq:principal} as
\begin{equation}\label{eq:rewrite-eq-principal-01}
\dfrac{\alpha^m}{4\sqrt{2}}-B_{n+1}^x=-B_n^x+\dfrac{\beta^m}{4\sqrt{2}}.
\end{equation}
Dividing both sides of equation  \eqref{eq:rewrite-eq-principal-01} by $B_{n+1}^{x}$, taking absolute value and using the inequality \eqref{eq:upper-bound_for-Bm/Bm+1}, we obtain
\begin{equation}\label{eq:upper-bound-Lambda}
\abs{\alpha^m (4\sqrt{2})^{-1} B_{n+1}^{-x}-1}< 2 \left(\dfrac{B_n}{B_{n+1}} \right)^x < \dfrac{2}{5.8^x}.
\end{equation}
Put
\begin{equation}\label{eq:def-Lambda-1}
\Lambda_1:=\alpha^m(4\sqrt{2})^{-1}B_{n+1}^{-x}-1.
\end{equation}
If $\Lambda_1=0$, we get $\alpha^m=4\sqrt{2}B_{n+1}^x$, so $\alpha^{2m}\in\ZZ$, which is false for all positive integers $m$, therefore $\Lambda_1\neq0$.

At this point, we will use Matveev's theorem to get a lower bound for $\Lambda_1$. We set $s:=3$ and we take
$$
\gamma_1:=\alpha,\quad \gamma_2:=4\sqrt{2},\quad \gamma_3:=B_{n+1}, \quad b_1:=m,\quad b_2:=-1, b_3:=-x.
$$
Note that $\gamma_1,\gamma_2,\gamma_3\in\QQ(\sqrt{2})$, so we can take $D:=2$. Since $h(\gamma_1)=(\log\alpha)/2$, $h(\gamma_2)=(\log32)/2$ and $h(\gamma_3)=\log B_{n+1}<n\log\alpha$, we can take $A_1:=\log\alpha$, $A_2:=\log32$ and $A_3:=2n\log\alpha$. Finally, inequality \eqref{eq:encad-m-terms-n-x} implies that $m>(n-2)x\geq x$, thus we can take $B:=m$. We also have  $B:=m\leq nx+2<(n+2)x.$ Hence, Matveev's theorem implies that
\begin{eqnarray}\label{eq:lower-bound-Lambda-1}
\log\abs{\Lambda_1} &\geq& -1.4\times 30^6\times 3^{4.5} \times 2^2 \times (1+\log 2)(\log\alpha)(\log32)(2n\log\alpha)(1+\log m)\nonumber\\
 &\geq & -2.1\times 10^{13} n(1+\log m).
\end{eqnarray}
The inequalities \eqref{eq:upper-bound-Lambda}, \eqref{eq:def-Lambda-1} and \eqref{eq:lower-bound-Lambda-1}, give that
$$
\begin{array}{rcl}
x<1.2\times 10^{13}n(1+\log m)<2.1 \times 10^{13}n\log m,
\end{array}
$$
where we used the fact that $1+\log m<1.7\log m$ for all $m\geq 5$. Together with the fact that $m<(n+2)x,$ we get that
\begin{equation}\label{eq:upper-bound-x-terms-m-n}
x<2.1 \times 10^{13}n\log ((n+2)x).
\end{equation}

%%%%%%%%%%%%%%%%%%%%%
\subsection{Small values of $n$}
%%%%%%%%%%%%%%%%%%%%%
We next treat the cases when $n\in[2,37]$. In this case, 
$$
x<2.1\times 10^{13}n\log((n+2)x)<7.8\times 10^{14} \log(46x)
$$
so $x<4 \times 10^{16}$. 

We next take another look at $\Lambda_1$ given by expression \eqref{eq:def-Lambda-1}. Put
$$
\Gamma_1:=m\log\alpha-\log(4\sqrt{2})-x\log B_{n+1}.
$$
Thus, $\Lambda_1=e^{\Gamma_1}-1$. One sees that the right-hand side of \eqref{eq:rewrite-eq-principal-01} is a number in the
interval $[-B_n^x,-B_n^x+1]$. In particular, $\Lambda_1$ is negative, which implies that $\Gamma_1$ is negative. Thus,
$$
0<-\Gamma_1<\dfrac{2}{5.8^x},
$$
so 
\begin{equation}\label{eq:Baker-1}
0<x\left(\dfrac{\log B_{n+1}}{\log \alpha }\right)-m+\left(\dfrac{\log(4\sqrt{2}}{\log \alpha}\right)<\dfrac{2}{5.8^x\log \alpha}.
\end{equation}
For us, inequality \eqref{eq:Baker-1} is 
$$
0<x\gamma-m+\mu<AB^{-x},
$$
where 
$$
\gamma:=\dfrac{\log B_{n+1}}{\log \alpha},\quad \mu=\dfrac{\log(4\sqrt{2})}{\log \alpha},\quad A=\dfrac{2}{\log \alpha},\quad B=5.8.
$$
We take $M:=4\times 10^{16}$.

The program was developed in PARI/GP running with $200$ digits. For the computations, if the first convergent such that $q > 6M$ does not satisfy the condition $\varepsilon > 0$, then we use the next convergent until we find the one that satisfies the condition.  In one minute all the computations were done. In all cases, we obtained $x\leq 77$.
A computer search with Maple revealed in less than one minute that there are no solutions to the equation \eqref{eq:principal} in the range $n\in[3,37]$ and $x\in[3,77]$.

%%%%%%%%%%%%%%%%%%%%%%%%%%%%
\subsection{An upper bound on $x$ in terms of $n$}
%%%%%%%%%%%%%%%%%%%%%%%%%%%%

From now on, we assume that $n \geq 38$. Recall from the previous section that %\eqref{eq:upper-bound-x-terms-m-n} , we have
\begin{equation}\label{eq:upper-bound-x-terms-m-n}
x<2.1 \times 10^{13}n\log ((n+2)x).
\end{equation}

Next we give an upper bound on $x$ depending only on $n$. If 
\begin{equation}\label{eq-upper-bound-s-cas1}
x\leq n+2,
\end{equation}
then we are through. Otherwise, that is if $n +2 < x$, we then have
$$
x<2.1 \times 10^{13}n\log x^2=4.2 \times 10^{13}n\log x,
$$
which can be rewritten as
\begin{equation}
\dfrac{x}{\log x}< 4.2 \times 10^{13}n.
\end{equation}
Using the fact that, for all $A\geq 3$
$$
\dfrac{x}{\log x}<A \quad \text{yields} \quad x<2A\log A,
$$
and the fact that $ \log(4.2\times 10^{13}n) <10\log n$, holds for all $n\geq 38$, we get that
\begin{eqnarray}\label{eq-upper-bound-s-cas2}
x &<& 2 (4.2\times 10^{13} n) \log((4.2\times 10^{13}n) \\
  &<& 8.4 \times 10^{13} n ( 10\log n)\nonumber\\
  &<& 8.4 \times 10^{14} n \log n.\nonumber
\end{eqnarray}
From \eqref{eq-upper-bound-s-cas1} and \eqref{eq-upper-bound-s-cas2}, we conclude that the inequality
\begin{equation}\label{eq:upper-bound-x-terms-m}
x < 8.4 \times 10^{14} n \log n
\end{equation}
holds for any $n\geq 38.$

%%%%%%%%%%%%%%%%%%%%%%%%
\subsection{An absolute upper bound on $x$}
%%%%%%%%%%%%%%%%%%%%%%%%

Let us look at the element
$$
y:=\dfrac{x}{\alpha^{2n}}.
$$
The above inequality \eqref{eq:upper-bound-x-terms-m}, implies that
\begin{equation}\label{eq:upper-bound-y}
y< \dfrac{8.4 \times 10^{14} n \log n}{\alpha^{2n}}< \dfrac{1}{\alpha^n},
\end{equation}
where the last inequality holds for any $n\geq 23$. In
particular, $y < \alpha^{-38} < 10^{-31}$. We now write
$$
B_n^x=\dfrac{\alpha^{nx}}{32^{x/2}} \left(1-\dfrac{1}{\alpha^{2n}} \right)^x,
$$
and 
$$
B_{n+1}^x=\dfrac{\alpha^{(n+1)x}}{32^{x/2}} \left(1-\dfrac{1}{\alpha^{2(n+1)}} \right)^x.
$$
We have
$$
1< \left(1-\dfrac{1}{\alpha^{2n}}\right)<e^y<1+2y,
$$
because $y<10^{-31}$ is very small. The same
inequality holds if we replace $n$ by $n + 1$.
We now follow the argument from \cite{Luca-Oyono:2011} to get that
\begin{equation}\label{eq:upper-bound-modif-Lambda2}
\abs{\alpha^{m-(n+1)x}32^{(x-1)/2}-(1-\alpha^{-x})}< \dfrac{32^{x/2}}{\alpha^{m+(n+1)x}}+2y(1+\alpha^{-x})< \dfrac{1}{2\alpha^n}+\dfrac{396y}{197}<\dfrac{3}{\alpha^n},
\end{equation}
where we used the fact that $32^{x/2}/(\alpha^{(n+1)x})\leq (4\sqrt{2}/\alpha^{38})^x<1/2$, $m\geq (n-2)x\geq n$ and $\alpha^x\geq \alpha^3 >197$, as well as inequality \eqref{eq:upper-bound-y}. Hence,
we conclude that
\begin{equation}\label{eq:upper-bound-Lambda-2}
\abs{\alpha^{m-(n+1)x}32^{(x-1)/2}-1}<\dfrac{1}{\alpha^x}+\dfrac{3}{\alpha^n} \leq \dfrac{4}{\alpha^l},
\end{equation}
where $l:=\min\{n,x\}$. We now set
\begin{equation}\label{eq:def-Lambda-2}
\Lambda_2:=\alpha^{m-(n+1)x}32^{(x-1)/2}-1.
\end{equation}
and observe that $\Lambda_2\neq 0$. Indeed, for if $\Lambda_2=0$, then
$\alpha^{2((n+1)x-m)}=32^{x-1}\in\ZZ$. which is possible only when
$(n +1)x = m$. But if this were so, then we would get $0=\Lambda_2=32^{(x-1)/2}-1$, which leads to the conclusion that $x=1$, which is not possible. Hence, $\Lambda_2\neq 0$. Next, let us notice that since $x\geq 3$ and $m \geq 38$, we have that
\begin{equation}
\abs{\Lambda_2} \leq \dfrac{1}{\alpha^3}+\dfrac{1}{\alpha^{38}}<\dfrac{1}{2},
\end{equation}
so that $\alpha^{m-(n+1)x}32^{(x-1)/2}\in[1/2,1]$. In particular,
\begin{equation}
\label{eq:x}
(n+1)x-m < \dfrac{1}{\log\alpha}\left(\dfrac{(x-1)\log32}{2}+\log 2\right) < x \left(\dfrac{\log 32}{2\log\alpha}\right) <  x
\end{equation}
and
\begin{equation}
\label{eq:xx}
(n+1)x-m > \dfrac{1}{\log\alpha}\left(\dfrac{(x-1)\log32}{2}-\log 2\right) > 0.9x-1.4>0.
\end{equation}
We lower bound the left-hand side of inequality \eqref{eq:def-Lambda-2} using again Matveev's theorem. We take 
$$s:=2,\; \gamma_1:=\alpha,\; \gamma_2:=4\sqrt{2},\; b_1:=m-(n+1)x,\; b_2:=x-1,$$ 
$$D := 2,\; A_1 := \log\alpha,\; A_2 := \log32,\; \mbox{ and } B := x.$$ 
We thus get that
\begin{equation}\label{eq:lower-bound-Lambda-2}
\log\abs{\Lambda_2} > -1.4\times 30^5 \times 2^{4.5} \times 2^2 (1+\log2) (\log\alpha)(\log32)(1+\log x).
\end{equation}
The inequalities \eqref{eq:upper-bound-Lambda-2} and \eqref{eq:lower-bound-Lambda-2},  give 
$$
l < 4 \times 10^{10} \log x.
$$
Treating separately the case $l=x$ and the case $l=n$, following the argument in \cite{Luca-Oyono:2011} we have that the upper bound
\begin{equation}\label{eq:bounding-x}
x<7\times 10^{28}.
\end{equation}
always holds.

%%%%%%%%%%%%%%%%%%%%%%%
\subsection{Reducing the bound on $x$}
%%%%%%%%%%%%%%%%%%%%%%%

Next, we take 
$$
\Gamma_2:=(x-1)\log(4\sqrt{2})-((n+1)x-m)\log\alpha.
$$ 
Observe that $\Lambda_2=e^{\Gamma_2}-1$, where $\Lambda_2$ is given by \eqref{eq:def-Lambda-2}. Since $\abs{\Lambda_2}<\dfrac{1}{2}$, we have that $e^{\abs{\Gamma_2}}<2$. Hence, 
$$
\abs{\Gamma_2}\leq e^{\abs{\Gamma_2}}\abs{e^{\Gamma_2}-1}<2\abs{\Lambda_2}<\dfrac{2}{\alpha^x}+\dfrac{6}{\alpha^n}.
$$
This leads to
\begin{equation}\label{eq:Legendre-1}
\abs{\dfrac{\log(4\sqrt{2})}{\log\alpha}-\dfrac{(n+1)x-m}{x-1}} <\dfrac{1}{(x-1)\log\alpha} \left(\dfrac{2}{\alpha^x}+\dfrac{6}{\alpha^n}\right).
\end{equation}
Assume next that $x > 100$. Then $\alpha^x > \alpha^{100} > 10^{33} > 10^4x$. Hence, we get that
\begin{equation}\label{eq:Legendre-2}
\dfrac{1}{(x-1)\log\alpha}\left(\dfrac{2}{\alpha^x}+\dfrac{6}{\alpha^n} \right) < \dfrac{8}{x(x-1)10^4 \log\alpha} < \dfrac{1}{2200(x-1)^2}.
\end{equation}
Estimates \eqref{eq:Legendre-1} and \eqref{eq:Legendre-2} leads to
\begin{equation}\label{eq:Legendre-3}
\abs{\dfrac{\log(4\sqrt{2})}{\log\alpha}-\dfrac{(n+1)x-m}{x-1}} < \dfrac{1}{2200(x-1)^2}.
\end{equation}
By a criterion of Legendre, inequality \eqref{eq:Legendre-3} implies
that the rational number $((n+1)x-m)/(x-1)$ is a convergent to $\gamma:=\log(4\sqrt{2})/\log\alpha$. Let 
$$
[a_0,a_1,a_2,a_3,a_4,a_5,a_6,\ldots]=[0,1,57,1,234,2,1,\ldots]
$$
be the continued fraction of $\gamma$, and let $p_k/q_k$ be it's $k$th convergent. Assume that $((n+1)x-m)/(x-1)=p_k/q_k$ for some $k$. Then, $x-1=dq_k$ for some positive integer $d$, which in fact is the greatest commun divisor of $(n+1)x-m$ and $x-1$. We have the inequality
$$
q_{54}> 7\times 10^{28}>x-1.
$$
Thus, $k\in\{0,\ldots,53\}$. Furthermore, $a_k\leq 234$ for all $k=0,1,\ldots , 53$. From the known properties of the continued fraction, we have that
$$
\abs{\gamma-\dfrac{(n+1)x-m}{x-1}}=\abs{\gamma-\dfrac{p_k}{q_k}}> \dfrac{1}{(a_k+2)q_k^2} \geq \dfrac{d^2}{236(x-1)^2} \geq \dfrac{1}{236(x-1)^2},
$$
which contradicts inequality \eqref{eq:Legendre-3}. Hence, $x\leq 100$.
\subsection{The final step}
To finish, we go back to inequality \eqref{eq:upper-bound-modif-Lambda2} and rewrite it as
$$
\abs{\alpha^{m-(n+1)x}32^{(x-1)/2}(1-\alpha^{-x})^{-1}-1}<\dfrac{3}{\alpha^n(1-\alpha^{-x}}<\dfrac{4}{\alpha^n}.
$$
Recall that $x\in[3,100]$ and from inequalities \eqref{eq:x} and \eqref{eq:xx}, we have that
$$
0.9x-1.4<(n+1)x-m<x.
$$
Put $t:=(n+1)x-m$. We computed all the numbers $\abs{\alpha^{-t}32^{(x-1)/2} (1+\alpha^{-x})^{-1}-1}$ for all $x\in[3,100]$ and all $t\in\left[\lfloor 0.9x-1.4\rfloor, \lfloor x\rfloor\right].$ None of them ended up being zero and the smallest of these numbers is $>10^{-1}$. Thus,
$1/10<3/\alpha^n$, or $\alpha^n<30$, so $n\leq 3$ which is false.

%%%%%%%%%%%%%%%%%%%%%%%%%%%%%%%


\begin{thebibliography}{1}
%%%%%%%%%%%%%%%%%%%%%%%%%%%%%%

\bibitem{Baker-Davenport:1969} A.~Baker and H.~Davenport, \emph{The equations $3x^2 - 2= y^2$ and $8x^2 - 7= z^2$},
The Quarterly Journal of Mathematics, {\bf 20} (1) 129--137, 1969.

%\bibitem{Bravo-Luca:2013} J.~J. Bravo and F.~Luca, \emph{Coincidences in generalized Fibonacci recurrences}, Journal of Number Theory, {\bf 133} (6) 2121--2137, 2013.

\bibitem{Dey-Rout:2017} P. K. Dey and S. S. Rout, \emph{Diophantine equations concerning Balancing and Lucas balancing numbers}, Arch. Math. \textbf{108} (2017) 29-43.

%\bibitem{Chaves-Marques-Togbe:2012} A. P. Chaves, D. Marques and A. Togb\'{e}, \emph{On the sum of powers of terms of a linear recurrence sequence}, Bull Braz. Math. Soc. \textbf{43(3)} (2012), 397--406.

\bibitem{Dujella-Petho:1998}
A. Dujella and A. Peth\H{o}, \emph{A generalization of a theorem of Baker and Davenport}, Quart. J. Math. Oxford Ser. (2) \textbf{49} (1998), no. 195, 291--306.

\bibitem{Finkelstein:1965}
R. P.~Finkelstein, \emph{The house problem}, American Math. Monthly. \textbf{72} (1965), 1082--1088.

\bibitem{Luca-Oyono:2011}
F. Luca and R. Oyono, \emph{An exponential Diophantine equation related to powers of two consecutive Fibonacci numbers}, Proc. Japan Acad. Ser. A, \textbf{87} (2011), 45--50.

%\bibitem{Marques-Togbe:2010} D. Marques and A. Togb\'{e}, \emph{On the sum of powers of two consecutive Fibonacci numbers}, Proc. Japan Acad. Ser. A, \textbf{86} (2010), 174–-176.

\bibitem{Matveev:2000}
E. M. Matveev, \emph{An explicit lower bound for a homogeneous rational linear form in the logarithms of algebraic numbers, II}, Izv. Math. \textbf{64(6)} (2000) 1217--1269.

\bibitem{Ruiz-Luca:2000}
C.A.G. Ruiz and F. Luca, \emph{An exponential Diophantine equation related to the sum of powers of two consecutive $k$-generalized Fibonacci numbers}, Coll. Math., \textbf{137(2)} (2014) 171--188.

\bibitem{RFLT:2018} S. E. Rihane, B. Faye, F. Luca and A. Togb\'e, \emph{On the exponential Diophantine equation $P_n^x + P_{n+1}^x=P_m$} Preprint 2018.

%\bibitem{Shorey-Tijdeman:1960} T.N. Shorey and R. Tijdeman, \emph{Exponential Diophantine Equations}, Cambridge Tracts in Mathematics, Cambridge University Press, Cambridge, \textbf{87} (1986).

\end{thebibliography}
\end{document}